\documentclass[letterpaper,12pt,notitlepage]{report}
\usepackage{enumerate}
\usepackage{enumitem}
\input{tcilatex}
\usepackage{cite}
\usepackage{graphicx}
\usepackage{color}
\usepackage{float}
\usepackage{tabularx,ragged2e,booktabs,caption}

\begin{document}
 
\title{High Order Accurate Solution of Poisson's Equation in Infinite Domains for Smooth Functions}%

\author{Christopher R. Anderson\thanks{UCLA Department of Mathematics, Los Angeles, CA. 90095-1555}
 \thanks{This material is based upon work \author{}supported in part by the Army Research Office under ARO grant W911NF-15-1-0149.}}
\date{\today}
\maketitle

\begin{abstract}
In this paper a method is presented for evaluating the convolution of the Green's function for the Laplace operator with a specified function $\rho(\vec x)$ at all grid points in a rectangular domain $\Omega \subset {\mathrm R}^{d}$ ($d = 1,2,3$), i.e. a solution of Poisson's equation in an infinite domain. 4th and 6th order versions of the method achieve high accuracy when $\rho ( \vec x )$ possesses sufficiently many continuous derivatives. The method utilizes FFT's for computational efficiency, and has a computational cost that is $\rm O (N  \log N)$ where $\rm N$ is the total number of grid points in the rectangular domain.  
\end{abstract}

\pagestyle{myheadings} \thispagestyle{plain}
\renewcommand\thesection{\arabic{section}}
\renewcommand\thesubsection{\thesection.\arabic{subsection}}

\section{Introduction}

In this paper we present a procedure for evaluating the values of the integral

\begin{equation}
\phi( \vec x ) = \int_{{\mathrm R}^{d}} {\rm G}( | \vec x - \vec s | ) \, \rho( \vec s ) \,  d \vec s \qquad \vec x \in \Omega \label{GreensFcnIntegral}
\end{equation}

\noindent
at all points in a rectangular domain $\Omega \subset {\mathrm R}^{d}$ ($d = 1,2,3$) associated with a grid with uniform mesh width in each coordinate direction.    
${\rm G}(\vec x)$ is the free-space Greens' function for the Laplace operator given by 
$\displaystyle {\rm G}(x) = {1 \over 2} \, | x | $ , $\displaystyle {\rm G}( \vec x) =  {1 \over {2 \pi}} \log ( | \vec x | )$, and $\displaystyle {\rm G}(\vec x) =  - {1 \over {4 \pi}  | \vec x | }$ for $d = 1,2,3$ respectively. It is assumed that $\rho (\vec x) $ is a function that 
has support contained in $\Omega$. Since ${\rm G}( \vec x )$ is the Green's function, this task is equivalent to determining the 
values of $\phi ( \vec x )$ for $\vec x \in \Omega$ where $\phi$ satisfies 

\begin{equation}
\Delta \phi = \rho \qquad \vec x \in {\mathrm R}^{d} \label{InfDomainPoisson}
\end{equation}

\noindent
with specified behavior at infinity matching (\ref{GreensFcnIntegral}), i.e. a solution of Poisson's equation in an infinite domain. 

There are many different techniques to evaluate (\ref{GreensFcnIntegral}) or (\ref{InfDomainPoisson}) that achieve high efficiency and accuracy when applied to specific types of problems. The procedure in this paper is specifically tailored to the class of problems where $\rho ( \vec x) $ is continuously distributed throughout $\Omega$ and possesses many continuous derivatives and where $\phi ( \vec x) $ is to be evaluated at all points of a rectangular grid. This is a rather different class of problems for which fast multipole and related methods  are especially useful for, namely problems where $\rho ( \vec x) $ is a sum of delta functions, or smoothed delta functions, or when $\rho ( \vec x) $ is non-zero only on lower dimensional surfaces such as those that arise in the solution of integral equations (see \cite{BeatsonGreengard1997} for an introduction and references).

The general idea of the method is that associated with the family of methods 
\cite{James197771} \cite{McCorquodale2005}  
\cite{McCorquodale:2007:LCA}  \cite{Serafini2005481}  \cite{wang1999efficient}  where the solution values are constructed by combining the values of two finite domain solutions; a solution to the non-homogeneous Poisson's equation with homogeneous boundary conditions and a solution of a homogeneous Poisson equation with non-homogeneous boundary conditions -- a ``harmonic correction''. The specific details of the method are very similar to those in \cite{James197771} but with different component computational procedures developed to exploit the differentiablilty of $\rho ( \vec x) $ to obtain high accuracy. In particular, the solution of the non-homogeneous Poisson's equation is accomplished with spectral accuracy using discrete sine transforms and the harmonic correction component is obtained with either 4th or 6th order accuracy using a deferred correction technique. The computational cost of the method is ${\rm O} ( {\rm N}  \log {\rm N})$ where $\rm N$ is the total number of grid points.

It is expected that the accuracy and efficiency of the method presented here is similar to other methods that are tailored for smooth right hand sides and obtain efficiency through the use of FFT's such as the ``expanding domain type'' \cite{CANDERrecursiveDom14} \cite{ANDERSON2016194}  and those that utilize some form of a modified Greens' functions \cite{Hejlesen2013458} \cite{Qiang2010313}. This is not surprising, as all methods exploit FFT's to solve the same target problem, but there are differences in the methods, and in particular, how the accuracy and efficiency are influenced by the choice of numerical parameters and the difficulty of method implementation. In the method presented here, the computational domain required is a small fraction larger than the original domain and so has a smaller memory footprint than the expanding domain type methods, especially when high accuracy is required.  Also, since there is no explicit smoothing incorporated into the Green's function representation, spectral accuracy can be obtained for the non-homogeneous Poisson equation solution component. There is a formal limitation of sixth order accuracy due to the accuracy with which the harmonic component is constructed, but, as the computational results will show, the consequences of this limitation may not be observed until very fine discretization, if at all. 

The term ``spectral accuracy'' has been used, and will be used, in the following sections. This term is used to describe an approximation to operators that have errors with bounds that are ${\rm O}( h^p )$ where $h$ is the magnitude of the grid mesh size and $p$ is a function of the number of continuous derivatives that the argument to the operator possesses. For example, when the Trapezoidal method with uniform mesh spacing is used to approximate integrals of periodic functions, the Euler-MacLaurin error bound  is ${\rm O}( h^{2 m + 2})$ for functions that are $2m + 2$ times continuously differentiable \cite{atkinson1978introduction}. Also, estimates of the computational cost will generally be given in terms of, $\rm N$, the total number of computational grid points. In stating these estimates, it is assumed that the number of grid panels in each coordinate direction is ${\rm O} ( {\rm N}^{1 \over d})$, but this assumption is not required for the success of the method; the method has no restrictions on the aspect ratio of the mesh sizes associated with the computational grid.

\section{General Strategy}
The general strategy of the procedure is independent of the dimension and consists of using the equivalence of evaluating ({\ref{GreensFcnIntegral}) at all the grid points in $\Omega$  to that of determining the values of a solution to (\ref{InfDomainPoisson}) in $\Omega$. The values of the solution to (\ref{InfDomainPoisson}) in $\Omega$ are obtained as the sum of the solutions to two finite domain problems, $\phi^{\ast}$, a solution of Poisson's equation in $\Omega$ with homogeneous Dirichlet boundary values,

\begin{equation}
\Delta \phi^{\ast} = \rho \quad {\rm for} \quad \vec x \in \Omega \quad  {\rm and} \quad \phi^{\ast}(\vec x) = 0  \quad  {\rm for} \quad  \vec x \in \partial \Omega \label{PhiStarEqn}
\end{equation}

\noindent
and a harmonic correction, $\phi_H$, a solution of Laplace's equation with boundary values given by the evaluation of (\ref{GreensFcnIntegral}),
 
 \begin{equation}
\Delta \phi_H = 0 \quad {\rm for} \quad \vec x \in \Omega \quad  {\rm and} \quad \phi_H(\vec x) =  \int_{\Omega} {\rm G}( | \vec x - \vec y | ) \, \rho( \vec s ) \,  d \vec s   \quad  {\rm for} \quad  \vec x \in \partial \Omega \label{PhiHeqn}
\end{equation}

It is assumed that the support of $\rho$ is a distance greater than a fixed value $\sigma$ away from the boundary of $\Omega$, an assumption that is easily satisfied by increasing the size of the computational domain so that the condition is satisfied. It is this assumption that 
leads to a procedure that avoids the complications arising from the singular nature of the Green's function when evaluating the boundary values required to solve (\ref{PhiHeqn}) in two and three dimensions. 

For every dimension, discrete sine transforms are used to construct spectrally accurate approximations to $\phi^{\ast}$ determined by (\ref{PhiStarEqn}).  
If an FFT based sine transform is used, and the number of panels in each coordinate direction of the discretization satisfies the constraints of the FFT routines employed,
then the a solution to (\ref{PhiStarEqn}) requires ${\rm O}( {\rm N}\log {\rm N} )$ arithmetic operations where ${\rm N}$ is the total number of grid points in the domain. Typically, the FFT constraints on the number of panels is that it be a product of small primes, a constraint that can also be satisfied by expanding the computational domain a few grid panels in one or more of the coordinate directions. 

The construction of the harmonic correction requires two computations, the first being the approximate evaluation of the boundary values, and the second solving (\ref{PhiHeqn}). Under the assumption that the support of $\rho$ is a fixed distance away from the boundary, then the use of the Trapezoidal method to evaluate the integral (\ref{GreensFcnIntegral}) on the boundary of the computational domain will be spectrally accurate, i.e. the rate of convergence of the approximation will be determined by the  number of continuous derivatives of $\rho$. 
The computational cost of accumulating all the boundary values one boundary point a time is ${\rm O}( {\rm N}) $ in one dimension,
${\rm O}( {\rm N}^{3\over 2}) )$ in two dimensions, and ${\rm O}( {\rm N}^{4\over 3}) )$ in three dimensions. However, by re-arrangement of the required sums, and use of FFT based fast discrete convolutions, the computational cost of evaluating (\ref{GreensFcnIntegral}) at each boundary point can be reduced to ${\rm O}( {\rm N}\log {\rm N} )$ in two and three dimensions.   

Given boundary values, efficiently constructing solutions with high accuracy to (\ref{PhiHeqn}) in two and three dimensions presents significant difficulties (an exact solution is used in one dimension). The procedure employed here is based upon deferred correction  \cite{KellerAndPereya1977}. Approximate solutions to (\ref{PhiHeqn}) with either 4th or 6th order accuracy can be obtained, and both with ${\rm O}( {\rm N}\log {\rm N} )$ computational work.  

In the next three sections we provide the details for each of the component calculations of the method; obtaining solutions of Poisson equation for $\phi^{\ast}$, the efficient accumulation of boundary values for $\phi_H$, and the determination of $\phi_H$ as the solution to Laplace's equation. To both reduce and simplify the presentation, we discuss each of these details for two dimensional problems and indicate the modifications required to obtain methods for either one or three dimensional problems.

\section{Sine transform solution of Poisson's equation}

Given a continuous function $\rho ( \vec x )$ with support that is a fixed distance $\sigma$ from the 
boundaries of a rectangular region $\Omega = [a, \,b] \times [c, \, d]$ the two dimensional 
problem for $\phi^{\ast}$ consists of determining high order accurate approximations to  

\begin{equation}
\Delta \phi^{\ast} = \rho \quad {\rm for} \quad \vec x \in \Omega \quad  {\rm and} 
\quad \phi^{\ast}(\vec x) = 0  \quad  {\rm for} \quad  \vec x \in \partial \Omega \label{PhiStarEqn2D}
\end{equation}

\medskip
\noindent
at all grid points  $(x_i, \, y_j) = 
(c + i\,h_x, \, d + j \,h_y)$ for $i = 0 \ldots {\rm M}_x$ and $j = 0 \ldots {\rm M}_y$ where the mesh spacing in each direction is uniform and given by  
$\displaystyle h_x = {{(b-a)} \over {\rm M_x}}$ and $\displaystyle h_y = {{(d-c)} \over {\rm M_y}}$.

Since $\rho(x)$ has support in $\Omega$ and homogeneous boundary conditions are specified for $\phi^{\ast}$ on $\partial \, \Omega$, a spectrally accurate solution of 
(\ref{PhiStarEqn2D}) can be obtained using discrete sine transforms. Specifically, let $\beta_{k_1,k_2}$ be the coefficients of the sine series approximation to $\rho (\vec x)$ determined by the two dimensional discrete sine transform, 
\begin{equation}
\rho (\vec x) \approx \tilde \rho( \vec x )  = \sum_{k_1 = 1}^{\rm {\rm M}_x-1} \sum_{k_2 = 1}^{\rm {\rm M}_y -1}  \beta_{k_1, k_2} \sin \left( {{k_1 \pi (x - a)} \over {(b-a)} } \right)
\sin \left( {{k_2 \pi (y - c)} \over {(d-c)} } \right) \label{rhoSinSeries}
\end{equation}
\noindent
where 
\begin{equation}
\beta_{k_1, k_2} =  \left({2 \over {b-a}} \right) \left({2 \over {d-c}} \right)  
\sum_{i = 1}^{{\rm M}_x-1}
\sum_{j = 1}^{{\rm M}_y-1}
\rho( x_i , \,  y_j ) \, \sin \left( {{k_1 \pi (x_i - a)} \over {(b-a)} } \right) \sin \left( {{k_2 \pi (y_j - c)} \over {(d-c)} } \right)  \, h_x \, h_y
\end{equation}

\medskip
\noindent
A spectrally accurate approximation to $\phi^{\ast} (\vec x)$ is given by 

\begin{equation}
 \phi^{\ast} (\vec x) =
 \sum_{k_1 = 1}^{{\rm M}_x-1} 
 \sum_{k_2 = 1}^{{\rm M}_y-1}
 \alpha_{k_1, \, k_2} \sin \left( {{k_1 \pi (x - a)} \over {(b-a)} } \right) \sin \left( {{k_2 \pi (y - c)} \over {(d-c)} } \right) \label{phiStarSoln2D}
\end{equation}

\medskip
\noindent
where 

\begin{equation}
\alpha_{k_1, k_2} = -\beta_{k_1, \, k_2} \left({{k_1^2 \, \pi^2} \over {(b-a)^2} } + {{k_2^2 \, \pi^2} \over {(d-c)^2} } \right)^{\, -1} \, 
\end{equation}

\medskip 
\noindent
The required coefficients $\beta_{k_1, k_2}$ and the grid values $\phi^{\ast} (x_i, \, y_j)$ can be obtained 
using invocations of an FFT based discrete sine transform with appropriate data scaling. The computational cost when using 
FFT based transforms is ${\rm O}( {\rm N}\log {\rm N})$. 

The spectral accuracy of this approximation follows from the facts that $\tilde \rho (\vec x)$ given by 
(\ref{rhoSinSeries}) is a spectrally accurate approximation of $\rho (\vec x)$, $\Delta^{\,-1}$ is a bounded linear operator, 
and (\ref{phiStarSoln2D}) is an exact solution of Poisson's equation with right hand side $\tilde \rho (\vec x)$.

\section{Efficient accumulation of boundary values}

In order to construct the harmonic correction it is necessary to determine the boundary values of $\phi_H$. Since
$\rho ( \vec x )$ has support that is a fixed distance $\sigma$ from the 
boundaries of $\Omega = [a, \,b] \times [c, \, d]$, the integrand occurring in 

\begin{equation}
\phi(\vec x) = \int_{\Omega} {1 \over 2 \pi} \log(| \vec x - \vec s | ) \, \rho( \vec s ) \,   d \vec s   \label{GreensFcnIntegral2D}
\end{equation}
\noindent
is as differentiable as $\rho$ when $\vec x \in \partial \, \Omega$. In addition, the integrand vanishes at the boundaries
to an order equal to the order of differentiablilty of the integrand, and so, as a consequence of the Euler-MacLaurin error 
bound for the Trapezoidal method the approximation 

\begin{equation}
\phi(\vec x_{\,i, \, j} ) \approx 
 \sum_{p = 1}^{{\rm M}_x-1} 
 \sum_{q = 1}^{{\rm M}_y-1}
  {1 \over 2 \pi} \log(| \vec x_{\,i,\,j} - \vec s_{p,q} | ) \, \rho( \vec s_{\,p,\,q} ) \,  h_x \, h_y \label{TrapOrig2D}
\end{equation}

\noindent
is a spectrally accurate. (Since the integrand vanishes at 
the computational boundaries, contributions to the approximation from boundary points are not included in (\ref{TrapOrig2D}).)

Accumulating the boundary values using (\ref{TrapOrig2D}) one boundary grid point at a 
time leads to a total computational cost of  ${\rm O}({\rm N}^{3 \over 2})$. Fortunately, this cost can be reduced to  
${\rm O}( {\rm N}\log {\rm N} )$ by exploiting the fact that the double sum can expressed as a single sum of the values of discrete convolutions --- the latter being evaluated efficiently using an FFT based discrete Fourier transform. For example, for grid points along the left edge of the domain, $(a,t_j), \, t_j = c + j\, h_y, \, j = 0 \ldots {\rm M}_y$, (\ref{TrapOrig2D}) can be expressed as 
\begin{equation}
 \sum_{p = 0}^{{\rm M}_x} \left( 
 \sum_{q = 0}^{{\rm M}_y}
  {1 \over 2 \pi}
  \log\left(\sqrt{(a - u_p)^2 \, + \, (t_j - v_q)^2 } \right) \, \rho(u_p, \, v_q ) \,  \, h_y \right)  h_x   \label{Trap2D}
\end{equation}
\noindent
where $u_p = a + p \, h_x$ and $v_q = c + q \, h_y$.  

\noindent
The sum (\ref{Trap2D}) can be written a sum of vectors,
\begin{equation}
 \sum_{p = 0}^{{\rm M}_x} \left( \vec w_p \right)  \, h_x   \label{Trap2Dconvolution}
\end{equation}

\noindent
where for each value of the index $p$, $\vec w_p$ is the vector of values resulting from the discrete convolution of the two functions in the variable $s$,   
\begin{equation}
{1 \over 2 \pi} \log\left(\sqrt{(a - u_p)^2 \, + \, s^2 } \right) \quad {\rm and} \quad \rho( u_p, s)
\end{equation}
sampled at the grid points of a uniform grid with mesh spacing $h_y$, i.e.

\begin{equation}
(w_p)_j =  \sum_{q = 0}^{{\rm M}_y} \, \log(\sqrt{(a - u_p)^2 \, + \, (t_j - s_q)^2 }) \, \rho( u_p, s_q) \, h_y \label{convolution}
\end{equation}
\noindent
where $t_j = c + j\, h_y$ and $s_q = c + q \,h_y$.

The double sum in (\ref{Trap2D}) can thus be accumulated as a sum of vectors resulting from discrete convolutions, which, if computed
using FFT's has a total computational cost of ${\rm O}({\rm M}_x \, {\rm M_y} \log {\rm M}_y )$. The boundary values along the other edges of the 
computational domain can be obtained similarly. The resulting total computational cost of evaluating the boundary values using fast discrete convolutions
is thus
${\rm O}({\rm N} \log {\rm N}^{1 \over d}  ) = {\rm O}( {\rm N}\log {\rm N} )$.

In one dimension the evaluation of the boundary values require the evaluation of single sums with $\rm N$ terms, so alternative evaluation procedures are not needed. In three dimensions, if the evaluation of the boundary values is done one grid point at a time, the resulting computational cost 
${\rm O}({\rm N}^{4 \over 3}))$. As in two dimensions, the computational cost can be significantly reduced by recognizing that the values on each face of the computational domain can be accumulated as the sum of vectors resulting from two-dimensional discrete convolutions over slices of $\rho (\vec x)$ on planes
parallel to the face. If these convolutions are carried out using FFT's then the cost of evaluation in three dimensions is reduced to $ {\rm O}( {\rm N}\log {\rm N} )$ 
as well. 

\section{High order accurate solutions of Laplace's equation}

The harmonic correction, $\phi_{H}$, is determined as a solution of   
\begin{eqnarray}
\label{harmonicEQ}
\Delta  u &  = & 0 \quad   \vec x \in \Omega  \\ 
u( \vec x ) & =  &g(\vec x) \quad \vec x \in \partial \Omega \nonumber 
\end{eqnarray}

\medskip
\noindent
where $\Omega$ is a rectangular region and the discretization is based upon a uniform grid in each coordinate direction. In one dimension, an exact solution, i.e. a linear function, is used. A high order accurate solution in two and three dimensions is substantially more difficult. 
In this section we describe the computational procedure used to create 4th and 6th order accurate approximate solutions for the two dimensional problem; the procedure for three dimensional problems is essentially identical but with additional terms in the operators. For completeness, the specifications of the operators for three dimensions are given at the end of the section. 

Let $D^2_s$ be the standard three point second order difference operator approximation to the second derivative $\displaystyle u_{ss}$ using mesh spacing $h_s$ for each coordinate $s = x$ and  $y$. The stencil for this one-dimensional operator is given by 

\begin{eqnarray}
{ 1 \over {h_s^{\,2}}}
\left[
\begin{array}{rrr}
1  & -2  & 1 \\
\end{array}
\right] \nonumber
\end{eqnarray}

\noindent
The two dimensional operator $D^2_r D^2_s$ obtained by applying $D^2_r$ to the result of applying $D^2_s$ has the stencil
\begin{eqnarray}
{ 1 \over {h_r^{\,2} \, h_s^{\,2} }}
\left[
\begin{array}{rrr}
1  & -2  & 1 \\
-2  & \, \, 4  & -2 \\
1  &  -2  & 1 
\end{array}
\right] \nonumber
\end{eqnarray}
and is a second order approximation to $u_{ssrr}$. 
\noindent
The five point discrete Laplace operator in two dimensions is given by 

\begin{equation}
\Delta_h \, u \, = \, (D^2_x + D^2_y) \, u 
\end{equation}
and is a second order approximation to $\Delta u = u_{xx}+u_{yy}$.

\medskip
\noindent
The foundation for fourth order approximations is the following relation 

\begin{eqnarray}
\left[ \Delta_h \right. & + & {h_x^2 \over 12} \left( D^2_x D^2_y  \right) 
+ {h_y^2 \over 12} \left( D^2_x D^2_y  \right)   \left. \right] u \qquad \label{4thOrderRelation2D} \\
& = & \nonumber \\
f  & + & {h_x^2 \over 12} \left( D^2_x f \right)  
+ {h_y^2 \over 12} \left( D^2_y f \right)  + {\rm O}({|h|^4}) \nonumber
\end{eqnarray}
\noindent
that holds for sufficiently differentiable functions that satisfy Poisson's equation with Dirichlet boundary conditions, i.e. 
\begin{eqnarray}
\Delta  u &  = & f \quad   \vec x \in \Omega \\ 
u( \vec x ) & =  &g(\vec x) \quad \vec x \in \partial \Omega \nonumber \label{NonHomogeneous2D}
\end{eqnarray}
\noindent
Here $|h| = \max( {h_x, h_y} ) $.

\medskip
\noindent
While one can establish (\ref{4thOrderRelation2D}) using Taylor's theorem, 
the origin of this relation is based upon an adaptation of the technique of deferred 
correction \cite{KellerAndPereya1977} and begins with the asymptotic expansion for $\Delta_h$, i.e. 

\begin{eqnarray}
\Delta_h u & = & \Delta u +  {h_x^2 \over 12} u_{xxxx} + {h_y^2 \over 12} u_{yyyy}  + \label{DiscreteLaplaceExp2D} \\
& + &  {h_x^4 \over 360} u_{xxxxxx}
+ {h_y^4 \over 360} u_{yyyyyy} + {\rm O}({|h|^6}) \nonumber 
\end{eqnarray}

\noindent
If $\Delta u = f$, then 
\begin{eqnarray}
u_{xxxx} = f_{xx} - u_{xxyy}  &\,& u_{yyyy} = f_{yy} - u_{yyxx}  \label{Drelations2D} \\ 
\nonumber 
\end{eqnarray}

\noindent 
Second order approximations for $u_{xxxx}$ and $u_{yyyy}$ are obtained by using the second and fourth derivative approximations, $D^2_s$ and $D^2_r D^2_s$ respectively in the relations (\ref{Drelations2D}). Inserting these second order approximations of the fourth derivatives of $u$ into (\ref{DiscreteLaplaceExp2D}) and rearranging terms results in the relation (\ref{4thOrderRelation2D}). A direct application of deferred correction would consist of using these second order approximations to the fourth derivatives to form a modified right hand side for an equation based upon (\ref{DiscreteLaplaceExp2D}) in which only $\Delta_h$ appears on the left hand side. The relation (\ref{4thOrderRelation2D}) is used because it leads to a procedure with equivalent computational cost and possesses better numerical stability properties.

A fourth order approximation to (\ref{harmonicEQ}) can thus be obtained by solving the linear system of equations associated
with the operator equation

\begin{equation}
\left( \Delta_h  +  {h_x^2 \over 12} \left( D^2_x D^2_y \right) 
+ {h_y^2 \over 12} \left( D^2_x D^2_y  \right)  
 \right) u  \label{harmonic4thOrder2D}
 = 0
\end{equation}

\noindent
The stencil width of the difference operators in (\ref{harmonic4thOrder2D}) is one, 
so the solution of the corresponding linear system is equivalent to solving a linear 
system of equations associated with the operator equation 

\begin{equation}
\left( \Delta_h  +  {h_x^2 \over 12} \left( D^2_x D^2_y  \right) 
+ {h_y^2 \over 12} \left( D^2_x D^2_y \right)  
 \right)_0 u  = \tilde g \label{harmonic4thOrderNH2D}
\end{equation}
where the subscript 0 indicates difference operators with homogeneous boundary conditions and $\tilde g$ is a grid function that is non-zero only at grid points immediately adjacent to the boundary with values obtained by transferring the 
contributions of the boundary values in the difference operators to the right hand side of (\ref{harmonic4thOrder2D}). 

To obtain a solution of (\ref{harmonic4thOrderNH2D}) one can exploit the fact that eigenvectors of the difference operators in  (\ref{harmonic4thOrderNH2D}) are the discrete sine functions, i.e. the discrete grid functions $S_{k_1,k_2}$ obtained by sampling the function 
\begin{equation}
\sin( k_1 \pi {(x-a) \over {(b-a)}})\sin( k_2 \pi {(y-c) \over {(d-c)}})
\end{equation}
at the nodes of the interior points of the computational grid. Thus, by using an FFT based discrete sine transform it is 
possible to construct a solution to (\ref{harmonic4thOrderNH2D}) with ${\rm O}({\rm N}\log {\rm N})$ computational work where N is the total number of grid points in the computational domain. 

\bigskip

A sixth order approximate solution to (\ref{harmonicEQ}) is determined by carrying out  
an additional step of deferred correction. The procedure consists of first constructing, $u^{\,(1)}$, a fourth order solution of (\ref{harmonic4thOrderNH2D}), and then solving the linear system of equations associated with the operator equation

\begin{eqnarray}
\left[ \Delta_h \right. & + & {h_x^2 \over 12} \left( D^2_x D^2_y \right) 
+ {h_y^2 \over 12} \left( D^2_x D^2_y  \right)  
 \left. \right] u^{\,(2)} \qquad \label{6thOrderRelation2D} \\
& = & \nonumber 
\left( {h_x^4 \over 240} + {{h_x^2 \, h_y^2} \over {144}}  \right) D^2_x D^2_x D^2_y  \, u^{\,(1)}  +
\left( {h_y^4 \over 240} + {{h_x^2 \, h_y^2} \over {144}}  \right) D^2_y D^2_x D^2_y  \, u^{\,(1)} 
\end{eqnarray}
 
\noindent 
for a sixth order approximation, $u^{\,(2)}$, to (\ref{harmonicEQ}). The origin of (\ref{6thOrderRelation2D}) 
is the incorporation of second order difference approximations to the sixth order derivatives occurring in the leading terms of the asymptotic error expansion to the operator (\ref{harmonic4thOrder2D}). In the derivation of the form of this latter asymptotic expansion it is important to note that
explicit use is made of the fact that $u$ is harmonic and equality of cross partial derivatives.  

The difference operators on right hand side of (\ref{6thOrderRelation2D}), obtained by applying $D^2_x$ and $D^2_y$ to $D^2_x D^2_y$ have stencil width two,
so the right hand side values at grid points immediately adjacent to the boundary cannot be determined by a direct application of the difference operators, i.e. one needs values outside of the computational domain. Satisfactory results were obtained with right hand side values at points immediately adjacent to the boundary being determined by cubic polynomial extrapolation of the right hand side values of (\ref{6thOrderRelation2D})  at interior points not immediately adjacent to the boundary. 

Since the operator on the left hand side of (\ref{6thOrderRelation2D}) is identical to that associated with the fourth order approximation the same solution procedure can be employed, i.e.  the solution of the linear system of equations associated with (\ref{6thOrderRelation2D}) can be obtained using an  FFT based discrete sin transform with computational work that is 
${\rm O}({\rm N}\log {\rm N})$. 

In three dimensions the seven point discrete Laplace operator in three dimensions is given by 

\begin{equation}
\Delta_h \, u \, = \, (D^2_x + D^2_y + D^2_z) \, u
\end{equation}
and is a second order approximation to $\Delta u = u_{xx}+u_{yy} + u_{zz}$.

A fourth order approximation to (\ref{harmonicEQ}) in three dimensions is obtained by solving the linear system of equations associated
with the operator equation

\begin{equation}
\left[ \Delta_h  +  {h_x^2 \over 12} \left( D^2_x D^2_y + D^2_x D^2_z \right) 
+ {h_y^2 \over 12} \left( D^2_x D^2_y + D^2_y D^2_z \right)  
+ {h_z^2 \over 12} \left( D^2_x D^2_z + D^2_y D^2_z \right) \right] u  \label{harmonic4thOrder3D}
 = 0
\end{equation}

\noindent
which, like the corresponding equation in two dimension (\ref{harmonic4thOrder2D}) can be solved using discrete sine transforms with ${\rm O}({\rm N}\log {\rm N})$ computational work. 

A sixth order approximate solution to (\ref{harmonicEQ}) is determined by carrying out  
an additional step of deferred correction. The procedure consists of first constructing, $u^{\,(1)}$, 
a fourth order solution of (\ref{harmonic4thOrder3D}), and then solving the 
linear system of equations associated with the operator equation

\begin{eqnarray}
\left[ \Delta_h \right. & + & {h_x^2 \over 12} \left( D^2_x D^2_y + D^2_x D^2_z \right) 
+ {h_y^2 \over 12} \left( D^2_x D^2_y + D^2_y D^2_z \right)  
+ {h_z^2 \over 12} \left( D^2_x D^2_z + D^2_y D^2_z \right) \left. \right] u^{\,(2)} \qquad \label{6thOrderRelation3D} \\
& = & \nonumber 
\left( {h_x^4 \over 240} + {{h_x^2 \, h_y^2} \over {144}}  \right) D^2_x D^2_x D^2_y  \, u^{\,(1)}  +
\left( {h_y^4 \over 240} + {{h_x^2 \, h_y^2} \over {144}}  \right) D^2_y D^2_x D^2_y  \, u^{\,(1)} \\
& \, & \left( {h_x^4 \over 240} + {{h_x^2 \, h_z^2} \over {144}}  \right) D^2_x D^2_x D^2_z  \, u^{\,(1)}  +
\left( {h_z^4 \over 240} + {{h_x^2 \, h_z^2} \over {144}}  \right) D^2_z D^2_x D^2_z  \, u^{\,(1)} \nonumber \\
& \, & \left( {h_y^4 \over 240} + {{h_y^2 \, h_z^2} \over {144}}  \right) D^2_y D^2_y D^2_z  \, u^{\,(1)}  +
\left( {h_z^4 \over 240} + {{h_y^2 \, h_z^2} \over {144}}  \right) D^2_z D^2_y D^2_z  \, u^{\,(1)} \nonumber
\end{eqnarray}
 
\noindent 
Similar to the two dimensional case, cubic polynomial extrapolation were used to obtain values of the right hand side (\ref{6thOrderRelation3D}) at points on the boundary and adjacent to the boundary. Solution of the discrete linear system can be obtained with  ${\rm O}({\rm N}\log {\rm N})$ computational work.

\section{Computational results}

The computational results on a three dimensional problem were selected to demonstrate the accuracy of the procedure, the achievable rates of convergence, and the dependency of the convergence rates upon the differentiablilty of $\rho$. Results are also presented that demonstrate the variation in the computed values as the computational domain is expanded for a fixed $\rho$ and mesh size. The latter is of interest because it is one measure of how well the discrete approximations inherit the feature of functions determined by (\ref{GreensFcnIntegral}) or (\ref{InfDomainPoisson}) that the values obtained in a given region does depend upon the size of the domain containing the support of $\rho (\vec x)$. Lastly, results are presented demonstrating the reduction in computation time with the use of multi-threaded programming techniques. 

The test function $\rho( \vec x)$ used was a translate of a radial function of the form 

\begin{equation}
{\rm B}_{\epsilon,\,p} ( \vec x ) = 
\left\{
\begin{array}{ccc}
\gamma_{\epsilon,\,p}\left( 1 - \left| {{\vec x}\over {\epsilon}} \right|^{\,2} \right)^{\,p} & \qquad & | \vec x | \le \epsilon \\
0 & \qquad & | \vec x | > \epsilon
\end{array}
\right.
\label{Bepsilon}
\end{equation}

\noindent 
where  $\gamma_{\epsilon,\,p}$ is determined so that (\ref{Bepsilon}) has unit integral. These functions are $p-1$ times continuously differentiable and have the property that a solution can be constructed to 
(\ref{InfDomainPoisson}) with $\rho(\vec x)$ given by ${\rm B}_{\epsilon,\,p} ( \vec x )$ \cite{CANDERpolyMollifiers14}. All calculations were carried out using FFT implementations provided by the FFTW3 software package 
\cite{FFTW05}.

The first set of results was obtained using a  computational domain  $[-1,1] \times [-1,1] \times [-1,1]$ with
$\rho ( \vec x ) = {\rm B}_{\epsilon,\,p} ( \vec x - \vec a)$
with $\vec a = ( {{1}\over{\sqrt{31}}} , 0.2, 0.1)$  and $\epsilon = 0.4$. Calculations were performed using values of $p$ chosen so that the $\rho( \vec x)$ was 0, 2, 4, 6, and 8 times continuously differentiable. For each value of $p$, solutions were determined for a sequence of mesh sizes and estimates of the rates of convergence were constructed. Plots of the maximal relative pointwise errors in $\phi( \vec x)$ for different values of differentiablilty are shown in Figure 1 as the mesh size $h= h_x = h_y = h_z$ is varied from 0.01 to .1. This variation in the
mesh size corresponds to varying the number of panels in each coordinate direction from 20 to 200. In Table 1 estimated rates of convergence based upon linear regression of log-log date values are given. 
The rates 
of convergence when $\phi_H$ was obtained using a 4th order method were based upon values of $h$ between 0.01 and 0.02, i.e. data values plotted at the left edge of Figure 1. The rates of convergence for the 6th order method were based upon values of $h$ between 0.02 and 0.1, i.e. the data values near the right edge of Figure 1.

\begin{figure}[H]
\centering
\includegraphics[height=4.0in,width=4.0in]{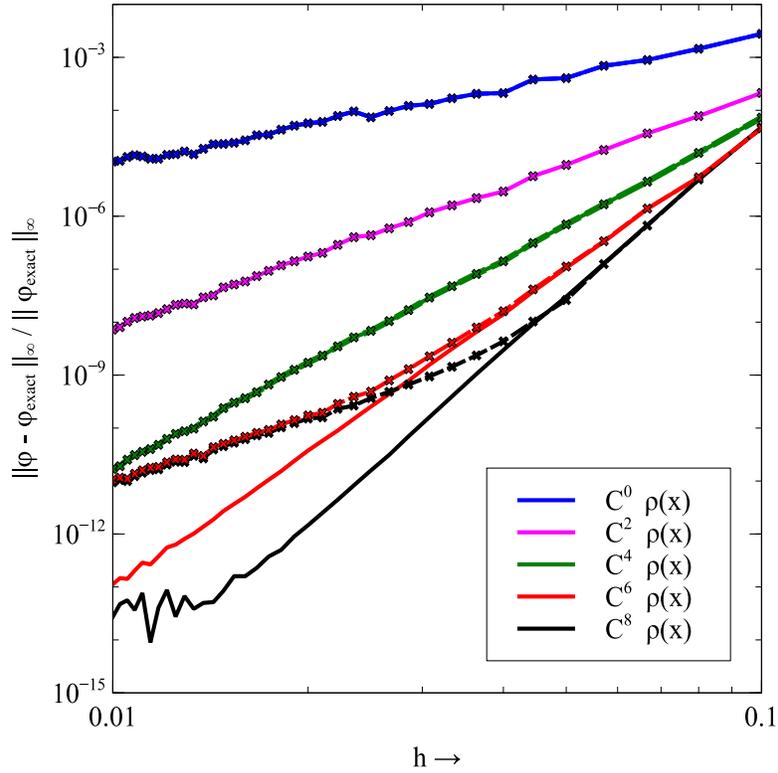} 
\caption{Convergence of approximate solution values of the 4th order (x-dashed lines) and 6th order (solid line) methods for $\rho(x)$ of differing degrees of differentiablilty.}
\end{figure}

It is apparent from these results that for this particular problem the errors associated 
with the computation of $\phi^{\ast}$ dominate for larger values of $h$; the method exhibits spectral accuracy that is expected of the method used to determine $\phi^{\ast}$. For smaller values of $h$ and when $\phi_H$ is computed using a 4th order method, the errors associated with the computation of $\phi_H$ dominate, and 4th order convergence is observed. When $\phi_H$ is computed using a 6th order method, the method still exhibits the spectral accuracy associated with $\phi^{\ast}$, indicating the errors are still dominated by errors in determining $\phi^{\ast}$. It is worth noting that when $\rho (\vec x)$ does not have many continuous derivatives, the pointwise error obtained for the 4th and 6th order method are essentially identical. Except for an increase in computational cost, there are no significant adverse consequences in using the sixth order method to determine the harmonic component for problems where $\rho ( \vec x) $ has limited differentiablilty.

{\renewcommand{\arraystretch}{1.2}
\begin{center}
\begin{tabular}{|l|c|c|c|c|c|}
\hline
$\rho$ \,\, differentiability &0&2&4&6&8 \\
\hline
${\rm{4th} \, \rm{order} \,\, \phi_H}$ &2.5&4.5&6.6&4.1&4.1 \\
\hline
${\rm{6th} \, \rm{order} \,\, \phi_H}$ &2.5&4.5&6.6&9.8&10.8 \\
\hline
\end{tabular}
\end{center}
\begin{center}
Table 1 : Estimated rates of convergence of $\phi$ as $h \to 0$.  
\end{center}

\bigskip
The second set of results were obtain by keeping the mesh width and $\rho(x)$ fixed and varying the size of the computational domain from 
an original $[-1,1] \times [-1,1] \times [-1,1]$ size to  $[-D,D] \times [-D,D] \times [-D,D]$ with $D \in [1,2]$. The function $\rho(x)$ was the same as that used to obtain the results shown in Figure 1 with $p$ chosen so that $\rho ( \vec x )$ was six times continuously differentiable. 
In Figure 2 are shown the relative maximal pointwise difference of values in $[-1,1] \times [-1,1] \times [-1,1]$ between a potential computed using the original computational domain and a potential computed using an extended computational domain.

\begin{figure}[H]
\centering
\includegraphics[height=4.0in,width=4.0in]{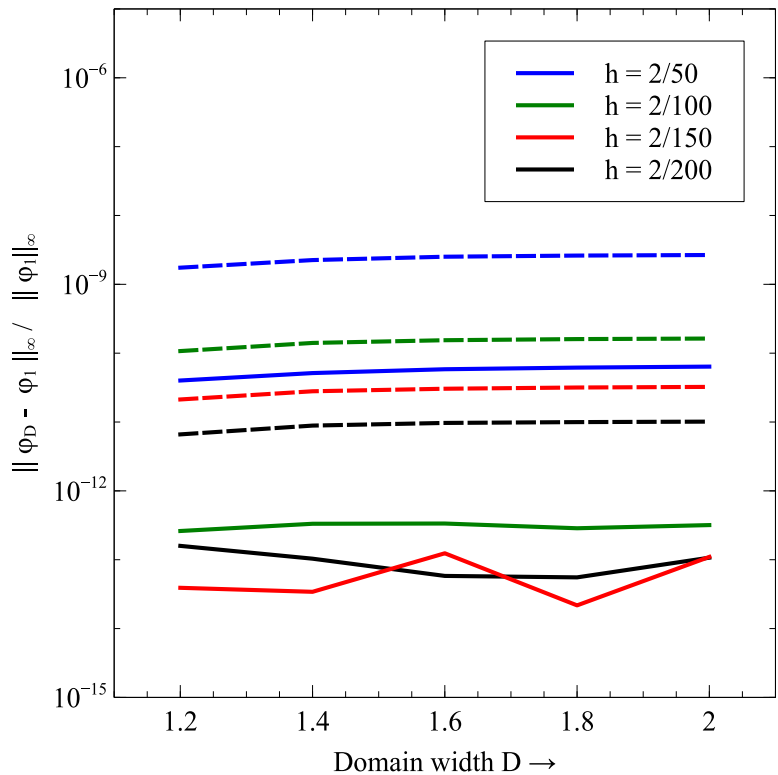} 
\caption{Maximal variation of values of potentials computed with increasing domain size. 4th order (dashed lines) and 6th order (solid line) harmonic correction.}
\end{figure}

The results presented in Figure 2 indicate that except for very small mesh sizes where numerical imprecision becomes significant, the fluctuations in the error as the domain size is increased are extremely small. Also the difference in solution values computed over the two different sized domains has a magnitude at most equal to the magnitude of the error associated with the approximation, and, with increasing 
mesh size becomes much smaller than the error associated with the approximation. In fact, except for the largest mesh size, the 6th order method exhibits,  up to essentially machine precision, the invariance of values with respect to domain size that is a feature of the exact solution.

\begin{figure}[H]
\centering
\includegraphics[height=4.0in,width=4.0in]{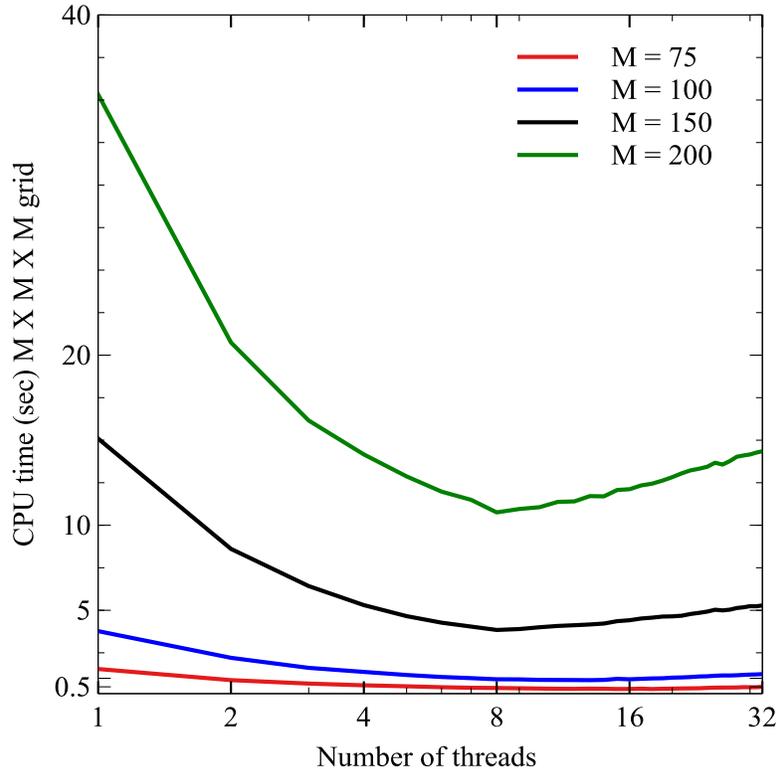} 
\caption{CPU time (seconds) as the number of threads is increased for problems of fixed size. }
\end{figure}

Many of computational tasks associated with the procedure can be decomposed into independent computational tasks that can be computed simultaneously, and thus it is possible to exploit multi-threaded programming on a multi-core CPU to improve computational performance. In particular, the most time consuming aspect of the procedure are the ${2 N_x} + {2 N_y} + {2 N_z} $ two dimensional convolutions required to accumulate the boundary values for the harmonic correction $\phi_H$. In the multi-threaded implementation whose results are shown in Figure 3 the loops which performed these computations were multi-threaded so that each thread carried out a separate convolution. The computational results shown in Figure 3 were obtained on a computer with a 32 core AMD Ryzen CPU and demonstrate the variation in computational efficiency as the number of threads is increased for problems of a fixed size. On all of the problems tested, a definite reduction in computational time is initially observed. The reduction is sub-linear, but with a more sophisticated multi-threading implementation it is likely that a nearly linear reduction in computational time could be obtained. Not unexpectedly, the performance deteriorates when using more than 8 threads. The important conclusion is that a multi-threaded implementation is certainly a way of reducing the computational time, but, as with all multi-threaded implementation, the optimal number number of threads to use is often less than the number of cores.

\section{Conclusion}

In this paper a procedure has been presented for constructing and evaluating at the grid points of a rectangular computational domain $\Omega$ a high order accurate approximation to the convolution of the  1, 2 and 3 dimensional Green's function of the Laplace operator with a smooth function $\rho ( \vec x)$ that has support contained within the interior of $\Omega$. Alternately, the procedure can be viewed as one that evaluates the restriction to $\Omega$ of  a high order accurate approximate solution of Poisson 's equation in an infinite domain.  The computational work required to obtain the solution is ${\rm O}({\rm N}\log {\rm N})$ where ${\rm N}$ is the total number of grid points in the computational domain. In one dimension the method is spectrally accurate, and in two and three dimensions versions of the procedure are formally 4th or 6th order accurate; although as the computational results show, one may observe substantially higher rates of convergence. Results on a three dimensional test problem demonstrate that for smooth functions the method can achieve the expected 4th and 6th order accuracy (or higher). For problems where the function $\rho( \vec x)$ is not smooth, the rate of convergence is reduced, but only to the rate corresponding to the accuracy of the spectral basis approximation of the solution of Poisson's equation in a finite domain.  Since many of the computational tasks involved in the construction of the procedure can be carried out independently implementations of the method can take advantage of multi-core CPU's to improve computational efficiency.  
\bigskip
\noindent

\medskip
\bibliographystyle{unsrt}
\bibliography{BibData}
\end{document}